\documentclass[10pt,twoside]{amsart}
\usepackage[english]{babel}
\usepackage[utf8]{inputenc}

\usepackage{amsmath, amscd, amssymb, mathrsfs, mathtools, fullpage, enumerate, svninfo,xcolor,enumitem}
\usepackage[all]{xy}
\usepackage{yfonts}
\usepackage{stmaryrd}
\usepackage{hyperref}

\usepackage{multicol}

\usepackage{tabu}

\usepackage{amsmath}
\usepackage{amssymb}

\numberwithin{equation}{section}

\newtheorem{thm}{Theorem}[section]

\newtheorem{prop}[thm]{Proposition}

\DeclareMathOperator{\TSym}{TSym}

\DeclareMathOperator{\Hom}{Hom}

\DeclareMathOperator{\GL}{GL}

\DeclareMathOperator{\Gal}{Gal}

\DeclareMathOperator{\Lie}{Lie}

\newcommand{\sH}{\mathscr{H}}

\newcommand{\sO}{\mathscr{O}}

\newcommand{\sP}{\mathscr{P}}

\newcommand{\cA}{\mathcal{A}}

\newcommand{\cD}{\mathcal{D}}

\newcommand{\cI}{\mathcal{I}}

\newcommand{\cO}{\mathcal{O}}
\newcommand{\cP}{\mathcal{P}}

\newcommand{\cS}{\mathcal{S}}

\newcommand{\fra}{\mathfrak{a}}
\newcommand{\frb}{\mathfrak{b}}
\newcommand{\frc}{\mathfrak{c}}

\newcommand{\frf}{\mathfrak{f}}

\newcommand{\frp}{\mathfrak{p}}

\newcommand{\CC}{\mathbb{C}}

\newcommand{\QQ}{\mathbb{Q}}

\newcommand{\ZZ}{\mathbb{Z}}

\newcommand{\NN}{\mathbb{N}}

\newcommand{\Q}{\mathbb{Q}}
\newcommand{\Z}{\mathbb{Z}}

\newcommand{\C}{\mathbb{C}}

\newcommand{\Zp}{{\ZZ_p}}

\newcommand{\Eis}{\mathrm{EK}}

\newcommand{\dR}{\mathrm{dR}}

\newcommand{\isom}{\cong}

\newcommand{\ol}[1]{\overline{#1}}
\DeclareMathOperator{\wP}{\widehat{\sP}}

\newcommand{\dashsum}{\sideset{}{'}{\sum}}

\newcommand{\EK}{EK}

\def\endpiece{xxx}
\def\makeAlphabet[#1]{\expandafter\makeA#1,xxx,}
\def\makealphabet[#1]{\expandafter\makea#1,xxx,}
\def\makeA#1,{\def\temp{#1}\ifx\temp\endpiece\else%
\mkbb{#1}\mkfrak{#1}\mkbf{#1}\mkcal{#1}\mkscr{#1}\expandafter\makeA\fi}%
\def\makea#1,{\def\temp{#1}\ifx\temp\endpiece\else\mkfrak{#1}\mkbf{#1}\expandafter\makea\fi}%
\def\mkbb#1{\expandafter\def\csname bb#1\endcsname{\mathbb{#1}}}%
\def\mkfrak#1{\expandafter\def\csname fr#1\endcsname{\mathfrak{#1}}}
\def\mkbf#1{\expandafter\def\csname b#1\endcsname{\mathbf{#1}}}
\def\mkcal#1{\expandafter\def\csname c#1\endcsname{\mathcal{#1}}}
\def\mkscr#1{\expandafter\def\csname s#1\endcsname{\mathscr{#1}}}
\def\makeop[#1]{\xmakeop#1,xxx,}
\def\mkop#1{\expandafter\def\csname #1\endcsname{{\mathrm{#1}}}} %
\def\xmakeop#1,{\def\temp{#1}\ifx\temp\endpiece\else\mkop{#1}\expandafter\xmakeop\fi}%

\makeop[Mic,Mod,im]

\newcommand{\OCp}{\sO_{\CC_p}}
\newcommand{\Cp}{{\CC_p}}
\newcommand{\Meas}{\mathrm{Meas}}

\AtBeginDocument{%
   \def\MR#1{}
}

\let\Re\relax
\DeclareMathOperator{\Re}{Re}

\begin{document}

\markboth{\hfill{\rm Guido Kings and Johannes Sprang} \hfill}{\hfill {\rm Algebraicity of critical Hecke $L$-values \hfill}}

\title{Algebraicity of critical Hecke $L$-values }

\author{Guido Kings and Johannes Sprang}

\begin{abstract}
In this survey, we review the known results on the algebraicity of critical values of Hecke $L$-functions and explain the new developments in \cite{Kings-Sprang}.
\end{abstract}

\maketitle

\setcounter{tocdepth}{1}

\section{Introduction}

 It is a remarkable fact that all special values of the Riemann zeta function $\zeta(s)=\sum_n n^{-s}$ at the even positive integers are given by the following formula:
\begin{equation*}
	\frac{\zeta(2n)}{(2\pi i)^{2n}}=-\frac{B_{2n}}{2(2n)!}\in \QQ,
\end{equation*}
where $B_n$ are the Bernoulli numbers defined by the generating series
\[
	\frac{z}{\exp(z)-1}=\sum_{n=0}^\infty B_n \frac{z^n}{n!}.
\]
It is an important insight by Deligne that the number $(2\pi i)^{2n}$ should be considered as a period of a motive attached to the zeta value $\zeta(2n)$. In fact, he conjectured a vast generalization of Euler's formula which predicts that all critical $L$-values divided by a certain period are algebraic. This conjecture of Deligne \cite{Deligne} is nowadays a special case of the vast Tamagawa number conjecture, which moreover predicts that the algebraic number can be expressed by arithmetic invariants of the motive. This more general conjecture is largely unknown for the motives attached to algebraic Hecke characters. 
In the following summary, based upon \cite{Kings-Sprang}, we discuss a complete solution of Deligne's conjecture in the case of Hecke $L$-functions, for Hecke characters which are attached to CM abelian varieties. It is conjectured that all Hecke characters can be obtained this way.  

For this, let $L$ be a number field of degree $d$ and $\frf\subseteq \sO_L$ an ideal. We write $J_L:=\Hom_{\QQ}(L,\CC)$ for the set of field embeddings of $L$ to $\CC$ and $\ZZ^{J_L}$ for the free abelian group on $J_L$. Let us write $\cI(\frf)$ for the group of fractional ideals of $L$ co-prime to $\frf$, and consider the following subgroup of principal ideals
\[
	\cP_{\frf}:=\{(\xi)\in \cI(\frf): \xi\equiv 1 \mod^\times \frf \text{ and } \sigma(\xi)>0 \text{ for all real embeddings }\sigma\in J_L \}.
\]
 An \emph{algebraic Hecke character} $\chi$ with values in a number field $E\subseteq \CC$ of conductor dividing $\frf$ and of infinity type $\mu=\sum_{\sigma \in J_L}\mu(\sigma)\sigma $ is a group homomorphism
\[
	\chi \colon \cI(\frf) \to E^\times,
\]
such that the restriction of $\chi$ to $\cP_\frf$ is given by the infinity type
\[
	\chi((\xi))=\xi^\mu:=\prod_\sigma \sigma(\xi)^{\mu(\sigma)}.
\]
The smallest ideal $\frf$ (with respect to divisibility) such that $\chi$ is an algebraic Hecke character modulo $\frf$ is called \emph{the  conductor of $\chi$}. The Hecke $L$-function of an Hecke character $\chi$ of conductor $\frf$ is
\[
	L(\chi,s):=\sum_{\substack{\fra\subseteq \sO_L \\ (\fra,\frf)=1}} \frac{\chi(\fra)}{N\fra^s}
\]
and the completed Hecke $L$-function is $\Lambda(\chi, s):=(|d_L|\cdot N\frf)^{s/2}\Gamma_\infty(\chi, s) \cdot L(\chi,s)$, where $\Gamma_\infty(\chi,s)$ is essentially a product of $\Gamma$-functions. The completed $L$-function satisfies a functional equation
\[
	\Lambda(\chi,s )=\epsilon(\chi)\Lambda(\chi^{-1},1-s)
\]
for some complex number $\epsilon(\chi)$ which has absolute value $1$ under each embedding. We call $L(\chi,s)$ \emph{critical} at an integer $s=k$ if and only if neither $\Gamma_\infty(\chi, s)$ nor $\Gamma_\infty(\chi^{-1}, 1-s)$ has a pole at $s=k$. By multiplication with the norm character $N(\cdot)$ and by the formula $L(\chi\cdot N(\cdot)^{-k},s)=L(\chi, s+k)$, the investigation of special values of Hecke $L$-functions is reduced to the study of these functions at $s=0$. We call the Hecke character $\chi$ critical if $L(\chi, s)$ is critical at $s=0$. Critical Hecke characters can only exist if the number field is either totally real or totally imaginary:
\begin{prop}
Let $\chi$ be a critical Hecke character of infinity type $\mu$ of a number field $L$, then either
\begin{enumerate}
\item[(i)] $L$ is totally real (i.e. all embeddings are real) and $\chi$ is of the form $\chi=\rho N(\cdot)^k$ for a finite Hecke character $\rho$ and some $k\in \ZZ$, or
\item[(ii)] $L$ is totally imaginary and contains a (maximal) CM subfield $K$. 
\end{enumerate}
\end{prop}
Recall that a \emph{CM field} $K$ is a quadratic totally imaginary extension of a totally real field. With this proposition, the discussion of critical Hecke $L$-values is divided into two cases. We remark right away that the new parts of our contribution concerns only the case $L\neq K$, as the totally real case is known by Siegel and Klingen \cite{Siegel}, \cite{Klingen} and the case of a CM field $L=K$ is due to Shimura \cite{Shimura} and Katz \cite{Katz-p-adic-L-fct}. The case of an imaginary quadratic field is already due to Damerell \cite{Damerell}. In the case of totally imaginary fields containing an imaginary quadratic field, algebraicity results have been obtained by Colmez and Bergeron--Charollois--Garcia, see \cite{Colmez,Bergeron-et-al}.

\section{Critical Hecke $L$-values}

We discuss the question about algebraicity of critical Hecke $L$-values up to explicit periods.

\subsection{Totally real fields}

The algebraicity of critical Hecke $L$-values for totally real fields is a celebrated theorem of Klingen and Siegel, which gives Deligne's conjecture in this case:

\begin{thm}[Siegel-Klingen, \cite{Siegel},\cite{Klingen}]
	Let $\chi=\rho N(\cdot)^{-k}$ be a critical Hecke character of a totally real field $L$, then
	\[
		\frac{L(\chi, 0)}{(2\pi i)^{\max(0,k)}}\in \ol{\QQ}.
	\]
	More precisely, the value is contained in the value field $E$ of $\chi$.
\end{thm}
We remark that more recently the ideas of equivariant Eisenstein cohomology classes, which are crucial in the approach in \cite{Kings-Sprang}, were successfully applied to reprove the Siegel-Klingen theorem, see \cite{BKL} and \cite{Bannai-et-al}.

\subsection{Totally imaginary fields}\label{subsec:tot-imag}

Let $\chi$ be a critical Hecke character of conductor $\frf$ and of infinity type $\mu$ of a totally imaginary field $L$. We define
\[
	\Sigma:=\{ \sigma\in J_L: \mu(\sigma)<0 \}.
\]
and write $\mu=\beta-\alpha$ with $\alpha\in \ZZ_{>0}^\Sigma$ and $\beta\in \ZZ_{\geq 0}^{\ol{\Sigma}}$.

One can show that $\Sigma$ allows us to view $\sO_L$ as a lattice of full rank in $\CC^\Sigma$ using the embedding:
\[
	\sO_L\hookrightarrow \CC^\Sigma, \quad \lambda\mapsto (\sigma(\lambda))_{\sigma \in \Sigma}.
\]
By CM theory, the associated complex torus $\CC^\Sigma/\sO_L$ are the $\C$-valued points of an abelian variety $\cA$ defined over
a number field $k$ and we assume that $k$ contains the Galois closure of $L$.

Our next goal is to recall the periods which are involved in the algebraicity statement: The action of $\sO_L$ on $\cA$ allows us to decompose the Lie algebra $\Lie(\cA/k)$ into eigenspaces corresponding to the different embeddings in $\Sigma$. We have
\begin{equation}\label{eq:decomp-Lie}
	\Lie(\cA/k)=\bigoplus_{\sigma\in \Sigma} \Lie(\cA/k)(\sigma),
\end{equation}
where $\Lie(\cA/k)(\sigma)$ is a one-dimensional $k$-vector space on which $\gamma\in \sO_L^\times$ acts by multiplication with $\sigma(\gamma)$. By passing to the dual, this decomposition induces a decomposition of the co-Lie algebra
\begin{equation}\label{eq:decomp-coLie}
	\omega_{\cA/k}=\bigoplus_{\sigma\in \Sigma} \omega_{\cA/k}(\sigma).
\end{equation}
For each $\sigma\in \Sigma$, let us fix a $k$-basis $\omega(\cA)(\sigma)\in  \omega_{\cA/k}(\sigma)$ of the $k$-vector space $\omega_{\cA/k}(\sigma)$. The Hodge filtration
\[
	\omega_{\cA/k}\subseteq H^1_{\dR}(\cA/k)
\]
allows us to view each $\omega(\cA)(\sigma)$ as an element of $H^1_{\dR}(\cA/k)$. 
\begin{def}\label{def:basis}
We write $\omega(\cA):=(\omega(\cA)(\sigma))_\sigma$ for the basis of $\omega_{\cA/k}$. 
\end{def}
Let us write
\[
	\langle\cdot, \cdot \rangle\colon H^1_{\dR}(\cA/k)\times H_1(\cA(\CC),\ZZ)\to \CC
\]
for the period pairing. Note that we have a canonical isomorphism $H_1(\cA(\CC),\ZZ)\cong \sO_L$, which gives us a distinguished element in $1\in H_1(\cA(\CC),\ZZ)$. The fixed basis of $\omega_{\cA/k}$ allows us to define periods
\[
	\Omega(\sigma):=\langle \omega(\sigma),1 \rangle\in \CC.
\]
We write 
\begin{equation}\label{eq:Omega}
    \Omega=(\Omega(\sigma))_{\sigma \in \Sigma}\in \CC^\Sigma.
\end{equation}
For $\alpha\in \NN^\Sigma$, let us define $\Omega^\alpha:=\prod_{\sigma\in \Sigma} \Omega(\sigma)^{\alpha(\sigma)}$. The basis of $\Lie(\cA/k)$ gives an isomorphism $\Lie(\cA/k)\otimes_k \CC\cong \CC^\Sigma$ and with this isomorphism the image of $H_1(\cA(\C),\Z)$ in $\CC^{\Sigma}$ becomes the period lattice $\Omega\sO_L$, so that
\begin{equation}\label{eq:period-lattice}
    \cA(\C)\isom \C^{\Sigma}/\Omega\sO_L.
\end{equation}

After choosing a basis of algebraic differential forms $\omega(\cA^\vee)(\sigma)$ of the dual abelian variety, one defines the periods $\Omega^\vee(\sigma)$, for details see \cite{Kings-Sprang}.

\begin{thm}[\cite{Kings-Sprang}]\label{thm:main}
Let $\chi$ be a critical Hecke character of $L$ of conductor $\frf$ and write its infinity type as $\mu=\beta-\alpha$ as above.
\begin{equation*}\label{eq:intro-L-value}
\frac{L(\chi,0)}{(2\pi i)^{-|\beta|}\Omega^{\alpha}\Omega^{\vee\beta}}\in \ol{\QQ}
\end{equation*}
More precisely, the value is contained in the compositum of the value field $E$ of $\chi$ and the field of definition of the abelian variety $\cA$ together with an auxiliary $\frf$-torsion point, for more details and much finer integrality results see \cite[Theorem 4.10]{Kings-Sprang}.
\end{thm}

The case where $L=K$ is a CM field has been previously obtained by Shimura \cite{Shimura} and Katz \cite{Katz-p-adic-L-fct}. The periods appearing in our main theorem, which are in this case the same as the periods of Shimura and Katz, are not the periods $c^+M(\chi)$ of the motive $M(\chi)$ attached to $\chi$ appearing in the Deligne conjecture. An important result of Blasius \cite{Blasius} relates these periods through the construction of a reflex motive and proves Deligne's conjecture in this case. 

In the general case $L\supset K$, Kufner \cite{Kufner} has extended the result of Blasius and deduced the full Deligne conjecture from our main result.

\begin{thm}[\cite{Kufner}]
Let $\chi$ be a critical Hecke character with value field $E$.
One has for the $E\otimes\CC$-valued $L$-value $(L(\tau\circ\chi,0))_{\tau\in J_E}\in E\otimes \CC$
\begin{equation*}
\frac{(L(\tau\circ\chi,0))_{\tau\in J_E}}{c^{+} R_{L/\Q}M(\chi)}\in E\subset E\otimes\CC,
\end{equation*}
where ${c^{+} R_{L/\Q}M(\chi)}$ is the period of the motive 
$R_{L/\Q}M(\chi)$ defined by Deligne in \cite{Deligne}.
\end{thm}

\section{$p$-adic interpolation}
A refined study of arithmetic invariants, in particular in Iwasawa theory, requires the study of congruences between special values of $L$-functions. These congruences are encoded in $p$-adic measures.

\subsection{$p$-adic measure}
Let $p$ be a prime, $G$ be a pro-finite group, $R$ a $p$-adically complete ring and $C(G,R)$ the continuous $R$-valued functions on $G$. The $R$-module of \emph{$p$-adic $R$-valued measure} is 
\[
	\Meas(G,R):=\Hom_{R}(C(G,R),R).
\] 
For $\mu \in \Meas(G,R)$ and $f\in C(G,R)$, we write
\[
	\int_Gf(g) d\mu(g) :=\mu(f)
\]
for the evaluation of $\mu$ on $f$. We equip $\Meas(G,R)$ with the convolution product: For $\mu,\nu\in \Meas(G,R)$, we define
\[
	(\mu\star \nu)(f):=\int_G\int_G f(g\cdot h) d\mu(g) d\nu(h).
\]
For $G=\Zp$ Mahler's description of $C(\Zp,R)$ implies a canonical isomorphism of $R$-algebras
\[
	\Meas(\Zp,R)\xrightarrow{\sim} H^0(\widehat{\mathbb{G}}_{m,R},\cO_{\widehat{\mathbb{G}}_{m,R}}),
\]
where $H^0(\widehat{\mathbb{G}}_{m,R},\cO_{\widehat{\mathbb{G}}_{m,R}})$ is the ring of formal functions on the formal multiplicative group $\widehat{\mathbb{G}}_{m,R}$ over $R$.

Let $\chi$ be an algebraic Hecke character of a number field $L$ of conductor $\frf$ co-prime to $p$. We fix an embedding $\ol{\QQ}\hookrightarrow \Cp$. To every algebraic Hecke character of conductor dividing $p^\infty \frf$ and infinity type $\mu\in \ZZ^{J_L}$, we can associate a $p$-adic Hecke character
\[
	\chi_p\colon \Gal(L(\frf p^\infty)/L)\to \CC_p^\times,
\]
as follows: Using the fixed embedding $\iota_p$, we consider $\chi$ as a character
\begin{equation}\label{eq:p-adic-Hecke-char-1}
	\chi \colon \cI(p\frf)\to \CC_p^\times,
\end{equation}
with
\begin{equation}\label{eq:p-adic-Hecke-char-2}
	\chi((\lambda))=\lambda^\mu,\quad \text{ for all }\lambda \in \cP_{p\frf}.
\end{equation}
Equation \eqref{eq:p-adic-Hecke-char-2} shows that $\chi(\cP_{p^n\frf})\subseteq 1+p^n\OCp$.
Together with the finiteness of $\cI(p\frf)/\cP_{p\frf}$, we deduce that $\chi(\cI(p\frf))\subseteq \OCp^\times$. Passing to the limit, we obtain a continuous homomorphism
\[
	\chi_p\colon \varprojlim_n \cI(p\frf)/\cP_{p^n\frf}\to \OCp^\times.
\]
By class field theory, the Artin map induces a canonical isomorphism
\[
	\varprojlim_n \cI(p\frf)/\cP_{p^n\frf}\xrightarrow{\sim} \Gal(L(p^\infty \frf)/L).
\]
The $p$-adic interpolation of critical $L$-values of Hecke characters consists of constructing a $p$-adic measure $\mu$ on the Galois group $\Gal(L(\frf p^\infty)/L)$ such that $\mu(\chi_p)$ is related to $L(\chi,0)$.

\subsection{Totally real fields}

The $p$-adic $L$-function of totally real fields has been constructed by Deligne--Ribet \cite{Deligne-Ribet}, Cassou-Nougès \cite{Cassou} and Barsky, \cite{Barsky}:
\begin{thm}
Let $L$ be a totally real field and $\frf \subseteq \sO_L$ an integral ideal. For every non-trivial auxiliary ideal $\frc \subseteq \sO_L$ there exists a $p$-adic measure $\mu_{\frf,\frc}\in \Meas(\Gal(L(\frf p^\infty)/L),\OCp)$ such that for every critical Hecke character of the form $\chi=\rho N(\cdot)^k$ of conductor dividing $\frf p^\infty$ with $k\geq 1$, we have
\[
	\int \chi_p(g) d\mu_{\frf,\frc}(g) = \{\text{explicit local factors}\}\cdot  L(\chi,0).
\]	
\end{thm}
More recent developments use equivariant cohomology classes for the construction of the $p$-adic $L$-function of totally real fields, see \cite{BKL} which is closely related to our approach in \cite{Kings-Sprang}. 

\subsection{Totally imaginary fields}

Already for imaginary quadratic fields there appears a distinction between split and inert primes in the construction of $p$-adic $L$-functions. For split primes, the corresponding $p$-adic measure has been constructed by Vishik-Manin \cite{Vishik-Manin} and Katz \cite{Katz-p-adic-Eisenstein}. In the inert case, the existence of a $p$-adic measure is not expected, however $p$-adic $L$-functions have been constructed in the weaker form of $p$-adic analytic functions or distributions. Let us, for example, mention Andreatta--Iovita, Boxall, Katz and Schneider--Teitelbaum, see \cite{Andreatta-Iovita}, \cite{Boxall}, \cite{Katz-formal}, \cite{Schneider-Teitelbaum}. Geometrically, the distinction between split and non-split primes corresponds to the reduction behaviour of the corresponding CM elliptic curve.

To make precise for which primes $p$ one can construct a $p$-adic measure interpolating the critical $L$-values for arbitrary totally imaginary fields, we need the following definitions:
Let $L$ be a totally imaginary field containing a (maximal) CM subfield $K$. A CM type of $K$ is a subset $\Sigma_K\subset J_K$ such that $J_K=\Sigma_K \coprod \overline{\Sigma}_K$, where $\overline{\Sigma}_K$ is the set of embeddings obtained from $\Sigma_K$ by complex conjugation. A CM type of $L$ is a subset $\Sigma\subset J_L$ of the form
\[
	\Sigma=\{\sigma\in J_L: \sigma|_K\in \Sigma_K\},
\]
for some CM type $\Sigma_K$ of the CM subfield $K$. We fix an embedding $\iota_p\colon\ol{\QQ}\hookrightarrow \Cp$. A CM type $\Sigma$ of $L$ is $p$-ordinary if the two sets of primes above $p$
\begin{align*}
	\Sigma_p:=\{\frp \text{ induced by the $p$-adic embeddings $\iota_p\circ \sigma$ with }\sigma\in\Sigma \}\\
	\overline{\Sigma}_p:=\{\overline{\frp} \text{ induced by the $p$-adic embeddings $\iota_p\circ \overline{\sigma}$ with }\overline{\sigma}\in\overline{\Sigma} \}.
\end{align*}
are disjoint. More explicitly, for a given prime $p$ there exists a $p$-ordinary CM type if and only if every prime above $p$ in $K^+$ splits completely in the extension $K/K^+$ where $K^+$ is the maximal totally real subfield. Such primes are called ordinary.

In the case $L=K$ a CM field, the corresponding $p$-adic measure for ordinary primes has been constructed by Katz, \cite{Katz-p-adic-L-fct}. For totally imaginary fields containing an imaginary quadratic field, $p$-adic $L$-functions have been constructed by Colmez--Schneps in \cite{Colmez-Schneps}.

\begin{thm}[\cite{Kings-Sprang}]\label{thm:main-p-adic}Let $\Sigma$ be a CM type of $L$, which is is ordinary for the prime number $p$.	For every fractional ideal $\frf$ there exists a $p$-adic measure $\mu_{\frf}$ on $\Gal(L(p^\infty\frf)/L)$ with the following interpolation property: For every critical algebraic Hecke character $\chi$ of attached CM type $\Sigma$ and  conductor dividing $p^\infty\frf$, we have
	\begin{equation*}
		\frac{1}{\Omega^{\chi}_p}\int_{\Gal(L(p^\infty\frf)/L)} \chi(g)d\mu_{\frf}(g)=\{\mbox{explicit local factors}\}\frac{L(\chi,0)}{\Omega^{\chi}}
	\end{equation*}
for an explicit $p$-adic period $\Omega_p^{\chi}$ which depends only on the infinity-type of $\chi$.
\end{thm}

\section{Methods} 

The main theorem of \cite{Kings-Sprang}, see Theorem \ref{thm:main}, generalizes earlier results of Shimura \cite{Shimura} and Katz \cite{Katz-p-adic-L-fct}, extending them from CM fields to arbitrary totally complex fields. We emphasize that our proof strategy differs from that of Katz even in the special case of a CM field $L = K$. Moreover, our method simultaneously produces all critical values of the $L$-functions, whereas Katz’ approach requires the use of the functional equation.

We now explain our approach to the critical values of Hecke $L$-functions and why the methods of Shimura and Katz do not extend directly. Both rely on real analytic Eisenstein series for the Hilbert modular group associated with the totally real subfield $K^+ \subset K$ and on the fact that their values at CM points can be related to $L$-values of Hecke characters. To establish algebraicity, they first show that holomorphic Eisenstein series are algebraic via the $q$-expansion principle. The real analytic Eisenstein series are then obtained by applying the Maaß–Shimura differential operators to the holomorphic ones. 
Katz had the important insight that the Maaß–Shimura operators arise naturally from the Gau\ss–Manin connection together with the Hodge decomposition. Their algebraicity at CM points follows then from the algebraicity of the Hodge splitting at such points. 
 
This approach breaks down for arbitrary extensions $L$ of $K$ of degree $n$, since the natural setting for the relevant Eisenstein series is the locally symmetric space attached to the arithmetic group $\GL_n(\sO_K)$, which is not Hermitian and therefore cannot be realized as the $\C$-points of an algebraic variety.

To overcome this difficulty, we have introduced several new perspectives in \cite{Kings-Sprang}. 
First, instead of working on the base of an abelian scheme, we work directly on the abelian scheme with CM. This idea goes back to Bannai and Kobayashi \cite{Bannai-Kobayashi}, who showed that Damerell’s result and Katz’ $p$-adic interpolation of $L$-values for CM elliptic curves can be recovered without the use of modular forms. Their key insight was that a certain theta function of the Poincaré bundle serves as a generating function for Eisenstein--Kronecker series at CM points. Our project began with the aim of understanding their method conceptually and extending it to higher-dimensional CM abelian varieties. This lead to the discovery of what we call 'Eisenstein--Kronecker classes'. Working on a single abelian variety means that $q$-expansions are no longer available to check algebraicity. 
A key feature of our approach is that algebraicity is built into the very definition of the Eisenstein–Kronecker classes. More precisely, we construct an equivariant coherent (and thus algebraic) cohomology class whose base change to $\C$ can be described by currents involving generalized Eisenstein–Kronecker series. 

A further obstacle is that a polarization on an abelian scheme with CM by $L$, where $L$ is a proper extension of a CM field $K$, is only linear over the maximal totally real subfield of $K$. We resolve this by working systematically with both the abelian scheme and its dual, which sheds new light on the periods appearing in formulas for special values of Hecke $L$-functions. This perspective also inspired Kufner’s new construction of Blasius’ reflex motive \cite{Kufner}. 

Another central insight of this paper is that one should not work with sections of algebraic bundles, as in the methods of Shimura, Katz, or Bannai–Kobayashi, but instead with equivariant coherent cohomology classes. This equivariant viewpoint is rooted in the ideas developed in \cite{BKL} and \cite{Graf}. 

Finally, our construction of $p$-adic $L$-functions differs substantially from Katz’ method in \cite{Katz-p-adic-L-fct}. In the CM case $L=K$, Katz investigates congruences of Eisenstein series through their $q$-expansions and deduces the existence of a $p$-adic measure. 
By contrast, our approach is inspired by the strategy of Bannai--Kobayashi. We use the equivariant cohomology class to construct a $p$-adic theta function of the Poincaré bundle. This $p$-adic theta function then yields the desired $p$-adic measure via Mahler’s theorem. 

%
%
 
\section{Sketch of the proof} 

In the following, we would like to sketch the proof of the Main Theorems of \cite{Kings-Sprang}. Let $\chi$ be a critical Hecke character of conductor $\frf$ and CM type $\mu$ of a totally imaginary field $L$. As explained in subsection \ref{subsec:tot-imag}, the Hecke character $\chi$ comes with a CM type $\Sigma$ of $L$ and we can write the infinity type as $\mu=\beta-\alpha$ with $\beta\in \ZZ_{\geq 0}^{\overline{\Sigma}}$ and $\alpha\in \ZZ_{>0}^{\Sigma}$. 

 A first step in the proof of Theorem \ref{thm:main} is to decompose the $L$-values into partial $L$-values
\[
	L(\chi,s)=\sum_{[\frb]\in \cI(\frf)/\cP_\frf} L(\chi,s,[\frb]),\quad \text{ with }L(\chi,s,[\frb]):=\sum_{\substack{\fra \in [\frb]\\ \fra \subseteq \sO_L}} \frac{\chi(\fra)}{N\fra^s},
\]
where $[\frb]$ runs over all classes in the ray class group $\cI(\frf)/\cP_\frf$. These partial $L$-values can then be related to Eisenstein series of the following form:
	\begin{equation}\label{def:Eisenstein}
		E^{\beta,\alpha}(t',s;\frf\frb^{-1},\sO_\frf^{\times}):=\dashsum_{\lambda\in \sO_\frf^\times \backslash (\frf\frb^{-1}+\sO_\frf^\times t')}\frac{\lambda^\beta}{{\lambda}^{\alpha}N(\lambda)^s},
	\end{equation}
	where $t'\in L$ and $\lambda$ runs over all non-zero $\sO_\frf^\times$-cosets of $\frf\frb^{-1}+\sO_\frf^\times t'$. Here, $\lambda^\alpha:=\prod_{\sigma\in \Sigma}\sigma(\lambda)^{\alpha(\sigma)}$ and $N(\lambda):=\prod_{\sigma\in \Hom(L,\CC)}\sigma(\lambda)$. The series \eqref{def:Eisenstein} converges a priori for $\Re(s)>[L:\QQ]-\frac{a}{2}$ with $a:=\sum_{\sigma\in \Sigma}\alpha(\sigma)$, but it can be extended to a meromorphic function on all of $\CC$, see \cite[\S 3.6]{Kings-Sprang}. In the case $L=K$ a imaginary quadratic field, the Eisenstein series \eqref{def:Eisenstein} are just the usual Eisenstein--Kronecker series. The partial $L$-function $L(\chi,s,[\frb])$ is now given in terms of the Eisenstein series by the formula
\begin{equation}\label{eq:L-fct-and-Eis}
    	L(\chi,s,[\frb])=\chi(\frb)N\frb^{-s} E^{\beta, \alpha}(1,s,\frf\frb^{-1},\sO_\frf^\times).
\end{equation}

To give these Eisenstein series a geometric interpretation, we work with the following setup. Let us keep the notation $\cA$ for algebraic abelian variety over some number field $k$ given by the CM type $\Sigma$, i.e. we have $\cA(\CC)=\CC^\Sigma/\sO_L$. The universal vectorial extension $\cA^\natural$ of the dual abelian variety $\cA$ parametrizes line bundles of degree $0$ with an integrable connection. The universal such line bundle $(\sP^\natural,\nabla)$ on $\cA\times \cA^\natural$ is called \emph{the Poincaré bundle with universal connection}.

 The main ingredient in the proof is an equivariant coherent cohomology class on an abelian scheme with values in $\wP^\natural$, the completion of the Poincar\'e bundle $\sP^\natural$ along $\cA\times e^{\natural}\subset \cA\times \cA^{\natural}$. Let $d=\dim \cA$ be the dimension of $\cA$ and 
$\Gamma\subset \sO_L^\times$  a subgroup of finite index which acts by automorphisms on the abelian scheme $\cA$. Then for a finite subscheme $\cD\subset \cA$ consisting of torsion sections and a $\Gamma$-invariant function $f$ on $\cD$ we construct a class
\begin{equation*}
\EK_{\Gamma}(f)\in H^{d-1}(\cA\smallsetminus\cD,\Gamma;\wP^\natural\otimes \Omega^{d}_{\cA/\cS} ),
\end{equation*}
which we call the \emph{equivariant coherent Eisenstein-Kronecker class}. 
It is absolutely essential for our applications to special values of Hecke $L$-functions to  have classes in equivariant cohomology. Otherwise, all specializations of this class along torsion sections would be trivial. The universal connection on the Poincaré bundle $(\sP^\natural,\nabla)$ allows us to take derivatives of $\EK_{\Gamma}(f)$ and after a suitable process of specialization along a torsion section $x\in (\cA\smallsetminus \cD)(k)$, we obtain classes:
\begin{equation*}
\Eis_{\Gamma}^{b,a}(f,x)\in H^{d-1}(\Gamma,H^{0}(k,\TSym^{a}(\omega_{\cA/k})\otimes\TSym^{b}(\sH)\otimes \omega^{d}_{\cA/k} )).
\end{equation*}
In a final step, we compute $\Eis_{\Gamma}^{b,a}(f,x)$ in terms of a $k$-basis of $\TSym^{a}(\omega_{\cA/k})\otimes\TSym^{b}(\sH)\otimes \omega^{d}_{\cA/k}$ (after taking the cap product with a generator of $H_{d-1}(\Gamma,\ZZ)$). One of these coefficients is
\[
	\frac{E^{\beta, \alpha}(1,s,\frf\frb^{-1},\sO_\frf^\times)}{(2\pi i)^{-|\beta|}\Omega^{\alpha}\Omega^{\vee\beta}}.
\]
which shows 
\[
	\frac{E^{\beta, \alpha}(1,s,\frf\frb^{-1},\sO_\frf^\times)}{(2\pi i)^{-|\beta|}\Omega^{\alpha}\Omega^{\vee\beta}}\in k.
\]
Let us briefly say a few words about the construction of the $p$-adic measure in Theorem \ref{thm:main-p-adic}. Let $p$ be a prime which is ordinary for the CM type $\Sigma$. Geometrically, this means that the corresponding abelian variety $\cA$ has good ordinary reduction at $p$. This implies that the formal scheme $\widehat{\cA}$, obtained by formal completion of the base change of $\cA$ to $\OCp$ at the zero section, is isomorphic to $\widehat{\mathbb{G}}_{m,\OCp}^d$. By restricting the equivariant cohomology class $\EK_{\Gamma}(f)\in H^{d-1}(\cA\smallsetminus\cD,\Gamma;\wP^\natural\otimes \Omega^{d}_{\cA/\cS} )$, to $\widehat{\cA}\times \widehat{\cA}^\vee$, we construct a $p$-adic theta function
\[
	\theta_\Gamma(f,x)\in H^0(\widehat{\cA}\times \widehat{\cA}^\vee, \cO_{\widehat{\cA}\times \widehat{\cA}^\vee})_\Gamma.
\]
Since $\widehat{\cA}\times \widehat{\cA}^\vee\cong \widehat{\mathbb{G}}_{m,\OCp}^{2d}$, we can use Mahler's isomorphism 
\[
	\Meas(\ZZ_p^{2d},\OCp)\xrightarrow{\sim} H^0(\widehat{\mathbb{G}}^{2d}_{m,\OCp},\cO_{\widehat{\mathbb{G}}^{2d}_{m,\OCp}}),
\]
to construct the $p$-adic measure.

\address{Fakult\"at f\"ur Mathematik, Universit\"at Regensburg, 93040 Regensburg, Germany\\
\email{guido.kings@mathematik.uni-regensburg.de}}

\address{Fakult\"at f\"ur Mathematik, Universit\"at Duisburg-Essen, 45127 Essen, Germany\\
\email{johannes.sprang@uni-due.de}}


\begin{thebibliography}{Yau}
\bibitem{Andreatta-Iovita}F. Andreatta and A. Iovi\c t\u a, Katz type $p$-adic $L$-functions for primes $p$ non-split in the ${CM}$ field, Comment. Math. Helv. {\bf 99} (2024), no.~4, 641--716; MR4815163

\bibitem{Barsky}D. Barsky, Fonctions zeta $p$-adiques d'une classe de rayon des corps de nombres totalement r\'eels, in {\it Groupe d'\'Etude d'Analyse Ultram\'etrique (5e ann\'ee: 1977/78)}, Exp. No. 16, 23 pp, Secr\'etariat Math., Paris, ; MR0525346


\bibitem{Bannai-Kobayashi}K. Bannai and S. Kobayashi, Algebraic theta functions and Eisenstein-Kronecker numbers, in {\it Proceedings of the Symposium on Algebraic Number Theory and Related Topics}, 63--77, RIMS K\^oky\^uroku Bessatsu, B4, Res. Inst. Math. Sci. (RIMS), Kyoto, ; MR2402003

\bibitem{Bannai-et-al}K. Bannai et al., Canonical equivariant cohomology classes generating zeta values of totally real fields, Trans. Amer. Math. Soc. Ser. B {\bf 10} (2023), 613--635; MR4583122


\bibitem{BKL}A.~A. Beilinson, G. Kings and A.~M. Levin, Topological polylogarithms and $p$-adic interpolation of $L$-values of totally real fields, Math. Ann. {\bf 371} (2018), no.~3-4, 1449--1495; MR3831278


\bibitem{Bergeron-et-al}N. Bergeron, P. Charollois and L.~E. Garcia, Transgressions of the Euler class and Eisenstein cohomology of ${\rm GL}_N({\bf Z})$, Jpn. J. Math. {\bf 15} (2020), no.~2, 311--379; MR4120422

\bibitem{Blasius}D. Blasius, On the critical values of Hecke $L$-series, Ann. of Math. (2) {\bf 124} (1986), no.~1, 23--63; MR0847951


\bibitem{Boxall}J.~L. Boxall, A new construction of ${\frp}$-adic $L$-functions attached to certain elliptic curves with complex multiplication, Ann. Inst. Fourier (Grenoble) {\bf 36} (1986), no.~4, 31--68; MR0867915

\bibitem{Cassou}P. Cassou-Nogu\`es, $p$-adic $L$-functions for totally real number field, in {\it Proceedings of the Conference on $p$-adic Analysis (Nijmegen, 1978)}, pp. 24--37, Katholieke Univ., Nijmegen, ; MR0522119

\bibitem{Colmez}P. Colmez, Alg\'ebricit\'e{} de valeurs sp\'eciales de fonctions $L$, Invent. Math. {\bf 95} (1989), no.~1, 161--205; MR0969417

\bibitem{Colmez-Schneps}P. Colmez and L. Schneps, $p$-adic interpolation of special values of Hecke $L$-functions, Compositio Math. {\bf 82} (1992), no.~2, 143--187; MR1157938


\bibitem{Damerell}R.~M. Damerell, $L$-functions of elliptic curves with complex multiplication. I, Acta Arith. {\bf 17} (1970), 287--301; MR0285540

\bibitem{Deligne}P. Deligne, Valeurs de fonctions $L$\ et p\'eriodes d'int\'egrales, in {\it Automorphic forms, representations and $L$-functions (Proc. Sympos. Pure Math., Oregon State Univ., Corvallis, Ore., 1977), Part 2}, pp. 313--346, Proc. Sympos. Pure Math., XXXIII, Amer. Math. Soc., Providence, RI, ; MR0546622

\bibitem{Deligne-Ribet}P. Deligne and K.~A. Ribet, Values of abelian $L$-functions at negative integers over totally real fields, Invent. Math. {\bf 59} (1980), no.~3, 227--286; MR0579702

\bibitem{Graf}P. Graf, Polylogarithms for $GL_2$ over totally real fields, arxiv.1604.04209


\bibitem{Katz-p-adic-Eisenstein}N.~M. Katz, $p$-adic interpolation of real analytic Eisenstein series, Ann. of Math. (2) {\bf 104} (1976), no.~3, 459--571; MR0506271

\bibitem{Katz-formal}N.~M. Katz, Formal groups and $p$-adic interpolation, in {\it Journ\'ees Arithm\'etiques de Caen (Univ. Caen, Caen, 1976)}, pp. 55--65, Ast\'erisque, No. 41--42, Soc. Math. France, Paris, ; MR0441928

\bibitem{Katz-p-adic-L-fct}N.~M. Katz, $p$-adic $L$-functions for CM fields, Invent. Math. {\bf 49} (1978), no.~3, 199--297; MR0513095

\bibitem{Katz-Another-look}N.~M. Katz, Another look at $p$-adic $L$-functions for totally real fields, Math. Ann. {\bf 255} (1981), no.~1, 33--43; MR0611271

\bibitem{Kings-Sprang}G. Kings and J. Sprang, Eisenstein--Kronecker classes, integrality of critical values of Hecke $L$-functions and $p$-adic interpolation, Ann. of Math. (2) {\bf 202} (2025), no.~1, 1--109; MR4927734


\bibitem{Klingen}H. Klingen, \"Uber die Werte der Dedekindschen Zetafunktion, Math. Ann. {\bf 145} (1961/62), 265--272; MR0133304


\bibitem{Kufner}H.~U. Kufner, Deligne's conjecture on the critical values of Hecke $L$-functions, arxiv.2406.06148


\bibitem{Nori}M.~V. Nori, Some Eisenstein cohomology classes for the integral unimodular group, in {\it Proceedings of the International Congress of Mathematicians, Vol.\ 1, 2 (Z\"urich, 1994)}, 690--696, Birkh\"auser, Basel, ; MR1403969


\bibitem{Schneider-Teitelbaum}P. Schneider and J.~T. Teitelbaum, $p$-adic Fourier theory, Doc. Math. {\bf 6} (2001), 447--481; MR1871671

\bibitem{Shimura}G. Shimura, On some arithmetic properties of modular forms of one and several variables, Ann. of Math. (2) {\bf 102} (1975), no.~3, 491--515; MR0491519

\bibitem{Siegel}C.~L. Siegel, \"Uber die Fourierschen Koeffizienten von Modulformen, Nachr. Akad. Wiss. G\"ottingen Math.-Phys. Kl. II {\bf 1970} (1970), 15--56; MR0285488

\bibitem{Vishik-Manin}M.~M. Vishik and Y.~I. Manin, $p$-adic Hecke series of imaginary quadratic fields, Mat. Sb. (N.S.) {\bf 95(137)} (1974), 357--383, 471; MR0371861



\end{thebibliography}
\end{document}